\documentstyle[12pt,amsmath,amssymb]{article}

\begin{document}
\newtheorem{lem}{Lemma}
\newtheorem{th}{Theorem}
\newtheorem{prop}{Proposition}
\newtheorem{rem}{Remark}
\newtheorem{define}{Definition}
\newtheorem{cor}{Corollary}

\allowdisplaybreaks

\newcommand{\N}{{\Bbb N}}
\newcommand{\C}{{\Bbb C}}
\newcommand{\Z}{{\Bbb Z}}
\newcommand{\R}{{\Bbb R}}
\newcommand{\Rp}{{\R_+}}
\newcommand{\FC}{{\cal F}C_{\mathrm b}(C_0(X),\Gamma_X)}
\newcommand{\FP}{{\cal FP}(C_0(X),\Gamma_X)}
\newcommand{\eps}{\varepsilon}
\newcommand{\la}{\langle}
\newcommand{\ra}{\rangle}
\newcommand{\D}{{\cal D}}
\newcommand{\fii}{\varphi}
\newcommand{\rom}[1]{{\rm #1}}
\newcommand{\dd}{\overset{{.}{.}}}

\begin{center}{\Large \bf
 The semigroup of the  Glauber dynamics of a continuous system of free particles
}\end{center}

{\large Yuri Kondratiev}\\
 Fakult\"at f\"ur Mathematik, Universit\"at
Bielefeld, Postfach 10 01 31, D-33501 Bielefeld, Germany;
Institute of Mathematics, Kiev, Ukraine; BiBoS, Univ.\ Bielefeld,
Germany.\\ e-mail:
\texttt{kondrat@mathematik.uni-bielefeld.de}\vspace{2mm}

{\large Eugene Lytvynov}\\ Department of Mathematics,
University of Wales Swansea, Singleton Park, Swansea SA2 8PP, U.K.\\
e-mail: \texttt{e.lytvynov@swansea.ac.uk}\vspace{2mm}

{\large Michael R\"ockner}\\
 Fakult\"at f\"ur Mathematik, Universit\"at
Bielefeld, Postfach 10 01 31, D-33501 Bielefeld, Germany;
 BiBoS, Univ.\ Bielefeld,
Germany.\\ e-mail: \texttt{roeckner@mathematik.uni-bielefeld.de }

{\small
\begin{center}
{\bf Abstract}
\end{center}

\noindent{We study properties of the  semigroup $(e^{-tH})_{t\ge
0}$ on the  space $L^ 2(\Gamma_X,\pi )$, where $\Gamma_X$ is the
configuration space over  a locally compact second countable
Hausdorff topological space $X$, $\pi$ is a Poisson measure on
$\Gamma_X$, and $H$ is the generator of the Glauber dynamics. We
explicitly construct the corresponding Markov semigroup of kernels
$({\bf P}_t)_{t\ge 0}$ and, using it, we prove the main results of
the paper: the  Feller property of the semigroup $({\bf
P}_t)_{t\ge 0}$ with respect to the vague topology on the
configuration space $\Gamma_X$, and the ergodic property of $({\bf
P}_t)_{t\ge 0}$. Following an idea of D.~Surgailis, we also give a
direct construction of the Glauber dynamics of a continuous
infinite system of free particles. The main point here is that
this process can start in every $\gamma\in\Gamma_X$, will never
leave $\Gamma_X$ and has cadlag sample paths in $\Gamma_X$.
}\vspace{1.5mm}

  \vspace{1.5mm}

 \vspace{1.5mm}

\section{Introduction}

The Glauber dynamics (GD) of a continuous infinite system of
particles, either free or interacting, is a special case of a
spatial birth and death process on the Euclidean space $\R^d$, or
on a more general topological space $X$. For a system of particles
in a bounded volume in $\R^d$, such processes were introduced and
studied by C.~Preston in \cite{P}, see also \cite{HS}. In the
latter case, the total number of particles is finite at any moment
of time.

 In the recent paper by L.~Bertini, N.~Cancrini, and F.~Cesi,
 \cite{BCC}, the generator of the GD  in
 a finite volume was studied. This generator corresponds to a special case of
 birth and death coefficients in Preston's dynamics. Under some
 conditions on the interaction between  particles, the authors
 of \cite{BCC} proved the existence of the spectral gap of the
 generator of the GD, which is uniformly positive with respect to
 all finite volumes $\Lambda$ and boundary conditions outside
 $\Lambda$. An explicit estimate of the spectral gap in a finite
 volume was derived by L.~Wu in \cite{Wu}.

The problem of construction of a spatial birth and death process
in the infinite volume was initiated by  R.A.~Holley and
D.W.~Stroock in  \cite{HS}, where it was solved in  a very special
case of nearest neighbor birth and death processes on the real
line.

In \cite{KL}, the GD in infinite volume was discussed. The process
now takes values in the configuration space $\Gamma_{\R^d}$ over
$\R^d$, i.e., in the space of all locally finite subsets in
$\R^d$, which is equipped with the vague topology. Using the
theory of Dirichlet forms \cite{MR1,MR2}, the authors of \cite{KL}
proved the existence of a Hunt process ${\bf M}$ on
$\Gamma_{\R^d}$ that is properly associated with the generator of
the GD with a quite  general pair potential of interaction between
particles. In particular, ${\bf M}$ is a conservative Markov
process on $\Gamma_{\R^d}$ with cadlag paths. An estimate of the
spectral gap of the generator of the GD in infinite volume was
also proved.

 In the case where the interaction between
particles is absent (i.e., the particles are free), the Poisson
measure $\pi$ on $\Gamma_{\R^d}$ is a stationary measure of the
GD. Let us recall that the Poisson measure possesses the chaos
decomposition property, and hence the space
$L^2(\Gamma_{\R^d},\pi)$ is unitarily isomorphic to the symmetric
Fock space over $L^2(\R^d)$, see e.g. \cite{S1}. It can be shown
that, under this isomorphism, the generator of the GD of free
particles goes over into the number operator $N$ on the Fock
space. The latter operator is evidently the  second quantization
of the identity operator, i.e., $N=d\operatorname{Exp}(\pmb 1)$.

On the other hand, a construction of a Markov process which
corresponds to the Poisson  space realization of the second
quantization of a doubly sub-Markov generator on $\R^d$ (or on a
more general  space) was proposed by D.~Surgailis in \cite{S2}.
However, D.~Surgailis did not discuss the following question: From
which configurations is  the process allowed to start so that it
never leaves the configuration space? It was only proved in
\cite{S2} that, for $\pi$-a.e.\ configuration $\gamma$, the
process starting at $\gamma$ will be a.s.\ in the configuration
space at some fixed time $t>0$.

In the case of the Brownian motion on the configuration space,
i.e., in the case of the independent motion of infinite Brownian
particles (cf.\ \cite{AKR3}), a solution to the above stated
problem was proposed by the authors  in \cite{KLR}. More exactly,
a subset $\Gamma_\infty$ of the configuration  space $\Gamma_X$
over a complete, connected, oriented, and stochastically complete
manifold $X$ of dimension $\ge2$ was constructed such that the
process can start at any $\gamma\in\Gamma_\infty$, will never
leave $\Gamma_\infty$, and has continuous sample paths in the
vague topology (and even in a stronger one). In the case of a
one-dimensional underlying manifold $X$, one cannot exclude
collisions of particles, so that a modification of the
construction of $\Gamma_\infty$ is necessary,  see \cite{KLR} for
details.

In this paper, we study properties of the semigroup of the GD of a
continuous infinite system of free particles. So, we fix a locally
compact second countable Hausdorff topological space $X$. We
denote by $\pi_m$ the Poisson measure on $\Gamma_X$ with intensity
$m$ being a Radon non-atomic measure on  $X$. In
Section~\ref{hbash}, we construct, on the space
$L^2(\Gamma_X,\pi_m)$, the Dirichlet form ${\cal E}$, the
generator $H$, and the semigroup $(e^{-tH})_{t\ge0}$ for the GD of
free particles in $X$. In particular, we derive an explicit
formula of the action of the semigroup on exponential functions
(Corollary~\ref{sdjsjdq}). The results of this section are
essentially  preparatory.

In Section~\ref{hsuer}, we construct a Markov semigroup of kernels
$({\bf P}_t)_{t\ge0}$ on $(\Gamma_X,{\cal B}(\Gamma_X))$ such
that, for each $F\in L^2(\Gamma_X,\pi_m)$, $$
(e^{-tH}F)(\gamma)=\int_{\Gamma_X} F(\xi)\, {\bf
P}_t(\gamma,d\xi),\qquad \text{$\pi_m$-a.e.\ }\gamma\in\Gamma_X$$
(Theorem~\ref{chbahs} and Proposition~\ref{cjj2356}). The reader
is advised to compare this result with \cite[Theorem~5.1]{KLR}.

The first  main result of the paper is Theorem~\ref{haajsd}, which
states that the semigroup $({\bf P}_t)_{t\ge0}$ possesses the
Feller property on $\Gamma_X$ (with respect to the vague
topology). Notice that, though we proved, in the case of the
Brownian motion
 on the configuration space, that the corresponding semigroup on
$\Gamma_\infty$ possesses a {\it modified\/} strong Feller
property \cite[Theorem~6.1]{KLR}, we were able to prove the usual
Feller property only in the case $X=\R^d$ and only with respect to
the intrinsic metric  of the Dirichlet form, which induces a
topology that is {\it much stronger\/} than the vague topology.

The second main result is that the semigroup $({\bf P}_t)_{t\ge
0}$ is ergodic (Theorem~\ref{uigvfuz}). More exactly, we show
that, for any probability measure $\mu$ on $\Gamma_X$, the image
measure ${\bf P}_t\mu$ of $\mu$ under ${\bf P}_t$ weakly converges
to $\pi_m$ as $t\to\infty$. Let us recall that, in the case of a
birth and death process in a finite volume in $\R^d$, the ergodic
property of a certain class of birth and death processes was
proved by C.~Preston \cite[Theorem~7.1]{P}. In particular, the
latter theorem holds in the case of the GD of free particles in a
finite volume.

Finally, in Section~\ref{ucdguz}, following the idea of
D.~Surgailis \cite{S2}, we give a direct construction of the
Glauber dynamics of an infinite system of free particles in $X$ as
a time homogeneous conservative strong Markov process on the state space $\Gamma_X$
with transition probability function $({\bf P}_t)_{t\ge0}$. The
main point here is that this process can start in every
$\gamma\in\Gamma_X$, will never leave $\Gamma_X$ and has cadlag
sample paths in $\Gamma_X$ (Theorem~\ref{tfut6t89}). In the case
where $m(X)=\infty$, we also show the possibility of restricting
the process to the subset $\Gamma_{X,\,{\mathrm inf}}$ of
$\Gamma_X$ consisting of all infinite configurations in $\Gamma_X$
(Corollary~\ref{wfdhci}).

\section{Generator of the Glauber dynamics}\label{hbash}

Let $X$ be a locally compact second countable Hausdorff
topological space. Such a space is known to be Polish, and we fix
a separable and complete metric $\rho$ on $X$ generating the
topology.  We denote by $C_0(X)$ the set of all continuous,
compactly supported, real-valued functions on $X$. Let ${\cal
B}(X)$ denote the Borel $\sigma$-algebra on $X$ and let $m$ be a
Radon measure on $(X,{\cal B}(X))$ without atoms.

We define  the configuration space $\Gamma_X$ over $X$ as the set
of all  locally finite  subsets of $X$:
$$\Gamma_X:=\{\gamma\subset X\mid |\gamma_\Lambda|<\infty\text{
for each compact }\Lambda\subset X\}.$$ Here, $|\cdot|$ denotes
the cardinality of a set and $\gamma_\Lambda:= \gamma\cap\Lambda$.
One can identify any $\gamma\in\Gamma_X$ with the positive Radon
measure $\sum_{x\in\gamma}\eps_x\in{\cal M}(X)$,  where  ${\cal
M}(X)$
 stands for the set of all positive
 Radon  measures  on
${\cal B}(X)$. We endow the space $\Gamma_X$  with the relative
topology as a subset of the space ${\cal M}(X)$ with the vague
topology, i.e., the weakest topology on $\Gamma_X$ with respect to
which  all maps
$$\Gamma_X\ni\gamma\mapsto\la\fii,\gamma\ra:=\int_X\fii(x)\,\gamma(dx)
=\sum_{x\in\gamma}\fii(x),\qquad\fii\in C_0(X),$$ are continuous.
We shall denote  the Borel $\sigma$-algebra on $\Gamma_X$ by
${\cal B}(\Gamma_X)$.

Let $\pi_{zm}$ denote the Poisson measure on $(\Gamma_X,{\cal
B}(\Gamma_X))$ with intensity $zm$, $z>0$. This measure can be
characterized by its Laplace transform
\begin{equation}\label{ewrwerewrwe5}\int_{\Gamma_X}
e^{\la\fii,\gamma\ra}\,\pi_{zm}(d\gamma)
=\exp\bigg(\int_X(e^{\fii(x)}-1)\,zm(dx)\bigg),\qquad \fii\in
C_0(X) .\end{equation} We refer e.g.\ to \cite{Kingman} for a
detailed discussion of the construction of the Poisson measure on
the configuration space. Let us recall that the Poisson measure
satisfies the Mecke identity:
\begin{equation}\int_{\Gamma_X} \pi_{zm}(d\gamma)\int_{X} \gamma(dx) \,
F(\gamma,x) =\int_{\Gamma_X} \mu(d\gamma)\int_{X} zm(dx)\,
F(\gamma\cup x,x)\label{fdrtsdrt}\end{equation} for any measurable
function $F:\Gamma_X\times X\to[0,+\infty]$, see \cite{Me67}. Here
and below, for simplicity of notation we just write $x$ instead
of $\{x\}$ for any $x\in X$.

It is easy to see that the Poisson measure $\pi_{zm}$ has all
local moments finite, i.e.,
\begin{equation}\label{sweqaw}\int_{\Gamma_X} \la \fii,\gamma\ra^n\,\pi_{zm}(d\gamma)<\infty,\qquad
\fii\in C_0(X),\ \fii\ge0,\ n\in\N.\end{equation}

We introduce the linear space  $\FC$ of all functions on $\Gamma_X$ of the
form
\begin{equation}\label{hbdv}
F(\gamma)=g_F(\la\fii_1,\gamma\ra,\dots,\la\fii_N,\gamma\ra),\end{equation}
where $N\in\N$, $\fii_1,\dots,\fii_N\in C_0(X)$, and $g_F\in
C_{\mathrm b}({\Bbb R}^N)$. Here, $C_{\mathrm
b}({\Bbb R}^N)$ denotes the set of all continuous bounded 
functions on $\R^N$.  For any $\gamma\in\Gamma_X$, we consider
$T_\gamma:=L^2(X,\gamma)$ as a ``tangent'' space to $\Gamma_X$
at the point $\gamma$, and for any $F\in\FC$ we define the
``gradient'' of $F$ at $\gamma$  as the element of $T_\gamma$
given by $D^-F(\gamma,x):=D^-_xF(\gamma):=F(\gamma\setminus
x)-F(\gamma)$, $x\in \gamma$. (Evidently, $D^-F(\gamma)$ indeed
belongs to $T_\gamma$.)

Let now $z:=1$.  We shall preserve the notation $\FC$ for the set
of all $\pi_m$-classes of functions from $\FC$. The set $\FC$ is
dense in $L^2(\Gamma_X,\pi_m)$. We  define
\begin{align} {\cal E}(F,G):=&
\int_{\Gamma_X}
(D^-F(\gamma),D^-G(\gamma))_{T_\gamma}\,\pi_m(d\gamma)\notag
\\=&\int_{\Gamma_X} \pi_m(d\gamma)\int_{X}\gamma(dx)\, D^-_x F(\gamma)\,
D^-_x G(\gamma), \qquad F,G\in\FC.\label{gcdtuf}\end{align} Notice
that, for any $F\in\FC$, there exists $\varphi\in C_0(X)$ such
that $|D^-_xF(\gamma)|\le \varphi(x)$ for all $\gamma\in\Gamma_X$
and $x\in\gamma$. Hence, due to  \eqref{sweqaw}, the right hand
side of \eqref{gcdtuf} is well defined. By \eqref{fdrtsdrt}, we
also get, for $F,G\in\FC$,
\begin{equation}\label{iud7u}{\cal
E}(F,G)=\int_{\Gamma_X}\pi_m(d\gamma)\int_{X} m(dx)\, D_x^+
F(\gamma)\, D_x^+ G(\gamma),\end{equation} where
$D_x^+F(\gamma):=F(\gamma\cup x)-F(\gamma)$.

Using \eqref{fdrtsdrt}, \eqref{sweqaw}, and \eqref{gcdtuf}, we see
that
\begin{equation}\label{hvgzfc} {\cal E}(F,G)=\int_{\Gamma_X}
(HF)(\gamma)G(\gamma)\, \pi_m(d\gamma),\qquad
F,G\in\FC,\end{equation} where
\begin{equation}\label{gen}
(HF)(\gamma)=-\int_{X}
D^+_xF(\gamma)\,m(dx)-\int_{X}D^-_xF(\gamma)\,\gamma(dx)
\end{equation}
and $HF\in L^2(\Gamma_X,\pi_m)$. Hence, the bilinear form $({\cal
E},\FC)$ is well-defined and closable on $L^2(\Gamma_X,\pi_m)$ and its closure will
be denoted by $({\cal E },D({\cal E}))$. Furthermore, using
  \cite[Theorem~4.1]{KL}, whose proof in the Poisson case admits a direct
generalization to the case of the general space $X$, we conclude
that the operator $(H,\FC)$ is essentially selfadjoint in
$L^2(\Gamma_X,\pi_m)$, and we denote  its closure by $(H,D(H))$. In
particular, $(H,D(H))$ is the generator of the bilinear form
$({\cal E },D({\cal E}))$.

Let us recall that the symmetric Fock space over $L^2(X,m)$ is
defined as the real Hilbert space given by  $$ {\cal
F}(L^2(X,m)):=\bigoplus_{n=0}^\infty {\cal F}_{n}(L^2(X,m)), $$
where ${\cal F}_0(L^2(X,m)):=\R$ and, for $n\in\N$, ${\cal
F}_n(L^2(X,m)):=L^2_{{\mathrm sym}}(X^n,m^{\otimes n})$  is the
subspace  of $L^2(X^n,m^{\otimes n})$ consisting of all symmetric
functions.

The Poisson measure $\pi_m$ possesses the chaos decomposition
property, and hence the space $L^2(\Gamma_X,\pi_m)$ is unitarily
isomorphic to the Fock space ${\cal F}(L^2(X,m))$, see e.g.\
\cite{S1}. More exactly, we set \begin{multline}\label{gxxdt}
{\cal F}(L^2(X,m))\ni f=(f^{(n)})_{n=0}^\infty \\ \mapsto
If:=\sum_{n=0}^\infty (n!)^{-1/2}
\int_{X^n}f^{(n)}(x_1,\dots,x_n)\, d{\cal X }(x_1)\dotsm d{\cal
X}(x_n)\in L^2(\Gamma_X,\pi_m).\end{multline} Here,
$\int_{X^n}f^{(n)}(x_1,\dots,x_n)\, d{\cal X }(x_1)\dotsm d{\cal
X}(x_n)$ denotes the $n$-fold multiple stochastic integral of the
function $f^{(n)}$ with respect to the Poisson random measure
${\cal X}(\Delta)(\gamma):=\gamma(\Delta)$, $\Delta\in{\cal
B}(X)$, $\gamma\in\Gamma_X$,  and the series on the right hand
side of \eqref{gxxdt} converges in $L^2(\Gamma_X,\pi_m)$. Then,
the operator $I$ is  unitary.

Denote  the number operator on ${\cal F}(L^2(X,m))$ by $N$, i.e.,
the domain of $N$ is given by $$
D(N)=\left\{f=(f^{(n)})_{n=0}^\infty:
\sum_{n=1}^\infty\|f^{(n)}\|_{{\cal F}_n(L^2(X,m))}^2 n^2<\infty
\right\}, $$ and $N\restriction {\cal F}_n(L^2(X,m))=n\pmb 1$,
$n\in\Z_+$, $\pmb 1$ being the identity operator. Evidently, $N$
is a positive selfadjoint operator  and, therefore, it generates
a contraction semigroup $(e^{-tN})_{t\ge0}$ in ${\cal
F}(L^2(X,m))$. In fact, $N$ is the
 second quantization of $\pmb 1$:
$N=d\operatorname{Exp}(\pmb 1)$, see e.g.\ \cite{RS,BK}, so 
$e^{-tN}=\operatorname{Exp}(e^{-t\pmb1})$.

\begin{prop}\label{zdfdf} Under the unitary operator  $I$\rom, the operator
$(N,D(N))$ goes over into $(H,D(H))$\rom.
\end{prop}

\noindent{\it Proof}. First, we note that, by the proof of
\cite[Lemma~3.2]{KL}, formula \eqref{iud7u} holds for all $F,G\in
D({\cal E})$. Denote by $\FP$ the set of all functions of the form
\eqref{hbdv}, where  $N\in\N$, $\fii_1,\dots,\fii_N\in C_0(X)$,
and $g_F$ is a polynomial on $\R^N$. By the proof of
\cite[Theorem~4.1]{KL}, $\FP\subset D(H)$, and furthermore, $H$ is
essentially selfadjoint on $\FP$. Analogously, we also  have that
$\FP\subset D(\cal E)$. Next, by
\cite[Theorem~5.1]{AKR3},\linebreak $\FP$ is a subset of the image
of $D(N)$ under $I$, and we have
\begin{align*}(INI^{-1}F,G)_{L^2(\Gamma_X,\pi_m)}&=\int_{\Gamma_X}\pi_m(d\gamma)
\int_X m(dx)\, D_x^+F(\gamma)D_x^+G(\gamma)\\& ={\cal E}(F,G)=
(HF,G)_{L^2(\Gamma_X,\pi_m)}\end{align*} for all $F,G\in\FP$.
Hence, $INI^{-1}$ coincides with $H$ on $\FP$. Since $INI^{-1}$ is
selfadjoint, this yields the statement. \quad $\square$

\begin{cor}\label{sdjsjdq}
Let $$ {\cal C}:=\{\varphi\in C_0(X)\mid -1<\varphi\le0\}.$$
Then\rom, for each $ \varphi\in{\cal C}$ and $t\ge0$\rom,
\begin{equation}\label{hssh}
e^{-tH}\exp[\la\log(1+\varphi),\cdot\ra]=\exp[\la\log(1+e^{-t}\varphi),\cdot\ra+(1-e^{-t})\la\varphi\ra]\quad
\text{\rom{$\pi_m$-a.e.}}\end{equation} Here\rom,
$\la\varphi\ra:=\int_X\varphi(x)\, m(dx)$\rom.
\end{cor}

\noindent{\it Proof}. The statement follows directly from
Proposition~\ref{zdfdf} by \cite[formula~(5.6)]{S1}. Let us
shortly repeat the arguments.

For each $\varphi\in C_0(X)$, we introduce the vector $$
\operatorname{Exp}\varphi:=\left(\frac1{n!}\,\varphi^{\otimes n
}\right)_{n=0}^\infty\in{\cal F}(L^2(X,m)),$$ which is called the
coherent state corresponding to the one-dimensional state
$\varphi$. We have, for any $\varphi\in{\cal C}$,
\begin{equation}\label{gzugfuzgf}I \operatorname{Exp}\varphi
=\exp[\la\log(1+\varphi),\cdot\ra-\la\varphi\ra]\quad
\text{$\pi_m$-a.e.},\end{equation} see e.g.\
\cite[formula~(5.5)]{S1}. Next, we evidently get
\begin{equation}\label{fdgdc}
e^{-tN}\operatorname{Exp}\varphi=\operatorname{Exp}(e^{-t}\varphi).\end{equation}
By \eqref{gzugfuzgf} and \eqref{fdgdc}, \eqref{hssh} follows from
Proposition~\ref{zdfdf}.\quad $\square$

\section{Semigroup of the Glauber dynamics}\label{hsuer}

We shall now construct probability  kernels for the semigroup
$(e^{-tH})_{t\ge0}$. Let $t>0$. Consider the two-point space
$\{0,1\}$  and define the probability measure $p_t$ on $\{0,1\}$
by setting $p_{t}(\{1\}):=e^{-t}$ and
$p_{t}(\{0\}):=1-e^{-t}$. Next,  let us fix an arbitrary
configuration $\gamma\in\Gamma_X$, $\gamma\ne\varnothing$. On the
product-space $\{0,1\}^\gamma$ equipped with the product
$\sigma$-algebra, we define the product-measure
$\bigotimes_{x\in\gamma}p_{t,x}$, where $p_{t,x}\equiv p_t$ for
each  $x\in\gamma$, . Set $$ \{0,1\}^\gamma\ni
a=(a(x))_{x\in\gamma}\mapsto {\cal U}(a)=\sum_{x\in\gamma}
a(x)\varepsilon_x\in\Gamma_X.$$
 It is easy to see that the mapping $\cal U$ is
measurable.  We denote by $P_{t,\gamma}$ the probability measure
on $\Gamma_X$ that is the image of the measure
$\bigotimes_{x\in\gamma}p_{t,x}$ under $\cal U$. We also set
$P_{t,\varnothing}:=\varepsilon_{\varnothing}$.

Now, for any $\gamma\in\Gamma_X$ and $t>0$, we define ${\bf
P}_{t,\gamma}$ as the probability measure on $\Gamma_X$ given by
the convolution of the measures $P_{t,\gamma}$ and
$\pi_{(1-e^{-t})m}$, i.e.,
\begin{equation}\label{gzgzuguguz}{\bf
P}_{t,\gamma}(A):=\int_{\Gamma_X}P_{t,\gamma}(d\eta_1)\int_{\Gamma_X}\pi_{(1-e^{-t})m}(d\eta_2)\,
\pmb 1_A(\eta_1+\eta_2),\qquad A\in{\cal
B}(\Gamma_X).\end{equation} Indeed, denote by $\dd\Gamma_X$ the
space of all $\N_0\cup\{\infty\}$-valued Radon measures on $X$ equipped with the
vague topology. Then, $\Gamma_X$ is a Borel-measurable subset of
$\dd\Gamma_X$. It is easy to see that the mapping
$$\Gamma_X\times\Gamma_X\ni(\eta_1,\eta_2)\mapsto
\eta_1+\eta_2\in\dd\Gamma_X$$ is measurable. Therefore,
considering $\pmb 1_A$ as an indicator function defined on
$\dd\Gamma_X$, we see that \eqref{gzgzuguguz} defines a measure on
$\Gamma_X$. To see that ${\bf P}_{t,\gamma}$ is a probability
measure, we have to show that, for $P_{t,\gamma}\otimes
\pi_{(1-e^{-t})m}$-a.e.\
$(\eta_1,\eta_2)\in\Gamma_X\times\Gamma_X$, $\eta_1+\eta_2$
belongs to $\Gamma_X$. By  construction, $P_{t,\gamma}$-a.e.\
$\eta_1\in\Gamma_X$ is a subset of $\gamma$. On the other hand,
the set $\gamma$ is of  zero $m$ measure and, therefore,
$\pi_{(1-e^{-t})m}$-a.e.\ configuration $\eta_2\in\Gamma_X$ has
empty intersection with $\gamma$, which implies the statement.

\begin{th}\label{chbahs} For each $F\in L^2(\Gamma_X,\pi_m)$ and $t>0$\rom, we
have \begin{equation}\label{kbhchj}
(e^{-tH}F)(\gamma)=\int_{\Gamma_X} F(\eta)\,{\bf
P}_{t,\gamma}(d\eta)\end{equation} for $\pi_m$-a\rom.e\rom.\
$\gamma\in\Gamma_X$\rom.
\end{th}

\noindent{\it Proof}. For $\gamma\in\Gamma_X$ and $\varphi\in{\cal
C}$, we have by \eqref{ewrwerewrwe5},  \eqref{gzgzuguguz}, and the
construction of $P_{t,\gamma}$:
\begin{align} &\int_{\Gamma_X}\exp[\la\log(1+\fii),\eta\ra]\,{\bf
P}_{t,\gamma}(d\eta)\notag\\ &\qquad
=\int_{\Gamma_X}P_{t,\gamma}(d \eta_1
)\int_{\Gamma_X}\pi_{(1-e^{-t})m}(d\eta_2)\exp[\la\log(1+\fii),\eta_1+\eta_2\ra]\notag\\
&\qquad=\int_{\Gamma_X} P_{t,
\gamma}(d\eta_1)\,\prod_{x\in\eta_1}(1+\fii(x))
\int_{\Gamma_{X}}\pi_{(1-e^{-t})m}(d\eta_2)
\exp[\la\log(1+\fii),\eta_2\ra] \notag\\ &\qquad
=\left(\prod_{x\in
\gamma\cap\operatorname{supp\varphi}}((1+\varphi(x))e^{-t}+(1-e^{-t}))
\right) \exp[(1-e^{-t})\la\varphi\ra]\notag\\
&\qquad=\exp[\la\log(1+e^{-t}
\fii),\gamma\ra+(1-e^{-t})\la\fii\ra],\label{hdfvbei}
\end{align}
where $\prod_{x\in\varnothing}c_x:=1$.

Next, for any measurable function $F:\Gamma_X\to[0,+\infty]$, we
have \begin{equation}\label{uacad} \int_{\Gamma_X}\int_{\Gamma_X}
 F(\eta)\,{\bf P}_{t,\gamma}(d\eta)\,\pi_m(d\gamma)=\int_{\Gamma_X}F(\gamma)\,
\pi_m(d\gamma).\end{equation} Indeed, it is easy to check that
$\{\exp[\la\log(1+\fii),\cdot\ra]\mid \fii\in{\cal C}\}$ is stable
under multiplication and  contains a countable subset separating
the points of $\Gamma_X$, so it generates ${\cal B }(\Gamma_X)$.
Therefore, we only have to check \eqref{uacad} for
$F=\exp[\la\log(1+\fii),\cdot\ra]$, $\varphi\in{\cal C}$. But this
immediately follows from \eqref{ewrwerewrwe5} and \eqref{hdfvbei}.

Now, by \eqref{hdfvbei} and \eqref{uacad}, the statement of the
theorem follows analogously to the proof of
\cite[Theorem~5.1]{KLR}. For the convenience of the reader, we repeat
 the arguments. It follows from \eqref{uacad} that, if
$A\in{\cal B}(\Gamma_X)$ is of zero $\pi_m$ measure, then ${\bf
P}_{t,\gamma}(A)=0$ for $\pi_m$-a.e.\ $\gamma\in\Gamma_X$.
Moreover, using the Cauchy--Schwarz inequality and \eqref{uacad},
we get \begin{align*}
\int_{\Gamma_X}\left(\int_{\Gamma_X}F(\eta)\,{\bf
P}_{t,\gamma}(d\eta)\right)^2\, \pi_m(d\gamma)&\le \int_{\Gamma_X}
\int_{\Gamma_X}|F(\eta)|^2\,{\bf
P}_{t,\gamma}(d\eta)\,\pi_m(d\gamma)\\
&=\int_{\Gamma_X}|F(\gamma)|^2\,\pi_m(d\gamma).\end{align*} Thus,
for each $t>0$, we can define a linear continuous operator $${\bf
P }_t:L^2(\Gamma_X,\pi_m)\to L^2(\Gamma_X,\pi_m)$$ by setting $$
({\bf P}_tF)(\gamma):=\int_{\Gamma_X}F(\eta)\, {\bf
P}_{t,\gamma}(d\eta).$$ By Corollary~\ref{sdjsjdq} and
\eqref{hdfvbei}, the action of the operator ${\bf P}_t$ coincides
with the  action of $e^{-tH}$ on the set
$\{\exp[\la\log(1+\fii),\cdot\ra]\mid \fii\in{\cal C}\}$, which is
total in $L^2(\Gamma_X,\pi_m)$ (i.e., its linear hull is a dense
set in $L^2(\Gamma_X,\pi_m)$). Hence, we get the equality ${\bf
P}_t=e^{-tH}$, which proves the theorem.\quad $\square$

In what follows, for a measurable function $F$ on $\Gamma_X$, we
set $$ ({\bf P}_tF)(\gamma):=\int_{\Gamma_X}F(\eta)\,{\bf
P}_{t,\gamma}(d\eta),\qquad t>0,\ \gamma\in\Gamma_X,$$ provided
the integral on the right hand side exists. Thus, by
Theorem~\ref{chbahs},  ${\bf P}_tF$ is a $\pi_m$-version of
$e^{-tH}$ for each $F\in L^2(\Gamma_X,\pi_m)$.

We also define a family of probability kernels $({\bf
P}_t)_{t\ge0}$ on the space $(\Gamma_X,{\cal B}(\Gamma_X))$
setting $$ {\bf P}_t(\gamma,A):={\bf P}_{t,\gamma}(A),\qquad
\gamma\in\Gamma_X,\ A\in{\cal B}(\Gamma_X),\ t\ge0.$$  Since
$\gamma\mapsto {\bf P}_tF(\gamma)$ is measurable for $F$ in the
linear span of $\{\exp[\la\log(1+\fii),\cdot\ra]\mid \fii\in{\cal
C}\}$ by \eqref{hdfvbei}, a monotone class argument shows that,
indeed, $\gamma\mapsto {\bf P}_{t}(\gamma,A)$ is ${\cal
B}(\Gamma_X)$-measurable for all $A\in{\cal B}(\Gamma_X)$.

\begin{prop}\label{cjj2356}
$({\bf P}_t)_{t\ge0}$ is a Markov semigroup of kernels on
$(\Gamma_X,{\cal B}(\Gamma_X))$\rom.
\end{prop}

\noindent{\it Proof}. We only have to prove the semigroup
property: ${\bf P}_t{\bf P}_s={\bf P}_{t+s}$, $t,s\ge0$. Let
${\varphi}\in{\cal C}$. By \eqref{ewrwerewrwe5} and
\eqref{hdfvbei}, we get, for any $\gamma\in\Gamma_X$,
\begin{align}& \int_{\Gamma_X}
\int_{\Gamma_X}\exp[\la\log(1+\varphi),\eta_1\ra]\, {\bf
P}_{s,\eta}(d\eta_1)\,{\bf P}_{t,\gamma}(d\eta)\notag\\ &\qquad=
\int_{\Gamma_X}\exp[\la\log(1+e^{-s}\fii),\eta\ra+(1-e^{-s})\la\varphi\ra]\,
{\bf P}_{t,\gamma}(d\eta) \notag\\ &\qquad=
\exp[\la\log(1+e^{-(t+s)}\fii),\gamma\ra+(1+e^{-t})e^{-s}\la\fii\ra+(1-e^{-s})\la\fii\ra]\notag\\
&\qquad= \exp[\la\log(1+e^{-(t+s)}\fii),\gamma\ra
+(1-e^{-(t+s)})\la\varphi\ra]\notag\\ &\qquad= \int_{\Gamma_X}
\exp[\la\log(1+\varphi),\eta\ra]\, {\bf
P}_{t+s,\gamma}(d\eta).\label{hscvv}
\end{align}
Analogously to the proof of \eqref{uacad}, we conclude from
\eqref{hscvv} that $$ \int_{\Gamma_X} \int_{\Gamma_X}F(\eta_1)\,
{\bf P}_{s,\eta}(d\eta_1)\,{\bf P}_{t,\gamma}(d\eta)
=\int_{\Gamma_X}F(\eta)\,{\bf P}_{t+s,\gamma  }(d\eta)$$ holds for
any measurable function $F:\Gamma_X\to[0,\infty]$, in particular,
for $F=\pmb 1_A$, $A\in{\cal B}(\Gamma_X)$.\quad $\square$

\begin{th}\label{haajsd} The semigroup $({\bf P}_t)_{t\ge0}$
possesses  the Feller property, i\rom.e\rom{.,} ${\bf
P}_t:C_{\mathrm b}(\Gamma_X)\to C_{\mathrm b}(\Gamma_X)$,
$t\ge0$\rom, where $C_{\mathrm b}(\Gamma_X)$ denotes the set of
all continuous bounded  functions on $\Gamma_X$\rom. \end{th}

\noindent{\it Proof}. First, we note that the space $\Gamma_X$ is
metrizable, see e.g. \cite{KMM}, and so the continuity of a
function $G:\Gamma_X\to\R$ follows if we can show the convergence
$G(\gamma_n)\to G(\gamma)$ as $\gamma_n\to\gamma$ vaguely in
$\Gamma_X$ .

Let us  fix any $x_0\in X$ and denote by $B(r)$ an open ball in
$X$ centered at $x_0$ and of radius $r>0$ with respect to the metric 
$\rho$. We recall that $\gamma_n\to\gamma$ vaguely in
$\Gamma_X$ if and only if, for any $r>0$, there exists $N\in\N$
such that $|\gamma_n\cap B(r)|=|\gamma\cap B(r)|{=:}l$ for all
$n\ge N$ and there exists a numeration of the points of the
configurations $\gamma_n\cap B(r)$, $n\ge N$, and $\gamma\cap
B(r)$, denoted by $\{x_k^{(n)}\}_{k=1}^l$, $\{x_k\}_{k=1}^l$,
respectively, such that $x_k^{(n)}\to x_k$ in $X$ as $n\to\infty$,
$k=1,\dots,l$. From here, by an  easy modification of the proof of
(6.9) in \cite{KLR}, we can conclude   the following

{\it Claim}. Let $\gamma_n\to\gamma$ as $n\to\infty$ vaguely in
$\Gamma_X$. Then, there exists a numeration of the points of the
configurations $\gamma_n$, $n\in\N$, and $\gamma$, denoted by $
\{x_k^{(n)}\}_{k\ge1}$, $\{x_k\}_{k\ge 1}$, respectively, such
that: 1) either for each $k\in\N$ if $|\gamma|=\infty$, or for
each $k\in\{1,\dots,|\gamma|\}$ if $|\gamma|<\infty$, there exists
$N\in\N$ such that $|\gamma_n|\ge k$ for all $n\ge N$ and
$x_k^{(n)}\to x_k$ in $X$ as $n\to\infty$.

So, we fix any $t>0$, $F\in C_{\mathrm b}(\Gamma_X)$ and
$\gamma_n$, $n\in\N$, $\gamma$ from $\Gamma_X $ such that
$\gamma_n\to\gamma$ vaguely, and we fix a numeration of the points
of $\gamma_n$, $n\in\N$, and $\gamma$ as described in the claim.

Analogously to the above, 
on the product-space
$\{0,1\}^\N$
we consider 
 the product measure $\bigotimes_{k\in\N} p_{t,k}$,
where $p_{t,k}\equiv p_t$. Then, as easily seen from
\eqref{gzgzuguguz},
\begin{equation}\label{shdo} ({\bf P}_t F)(\gamma_n)=\int_{\{0,1\}^\N}
\bigotimes_{k\in\N}
p_{t,k}(d(a_1,a_2,\dots))\int_{\Gamma_X}\pi_{(1-e^{-t})m}(d\eta)\,F
\left(\sum_{k=1}^{|\gamma_n|}a_k\varepsilon_{x_k^{(n)}}
+\eta\right).\end{equation} Set 
$D:=\left(\bigcup_{n=1}^\infty\gamma_n\right)\cup\gamma$. Since
$D$ has zero $m$ measure, $ \pi_{(1-e^{-t})m }$-a.e.\
$\eta\in\Gamma_X$ has empty intersection with $D$. Furthermore,
for any $\eta\in\Gamma_{X}$ with empty intersection with $D$ and
for any $(a_1,a_2,\dots)\in\{0,1\}^\N$, we easily get, using the
claim, that
\begin{equation}\label{gavusgv}
\sum_{k=1}^{|\gamma_n|}a_k\varepsilon_{x_k^{(n)}} +\eta\to
\sum_{k=1}^{|\gamma|}a_k\varepsilon_{x_k} +\eta\quad \text{vaguely
in $\Gamma_X$ as $n\to\infty$}.\end{equation} By the monotone
convergence theorem, we now have from \eqref{shdo} and
\eqref{gavusgv}: \begin{align*} ({\bf P}_t F)(\gamma_n)&\to
\int_{\{0,1\}^\N} \bigotimes_{k\in\N}
p_{t,k}(d(a_1,a_2,\dots))\int_{\Gamma_X}\pi_{(1-e^{-t})m}(d\eta)\,F
\left(\sum_{k=1}^{|\gamma|}a_k\varepsilon_{x_k}
+\eta\right)\\&\quad= ({\bf P}_t F)(\gamma)\quad\text{as
}n\to\infty,\end{align*} which yields the theorem.\quad $\square$

Before formulating the next theorem, let us recall that we use the
notation $\dd\Gamma_X$ for the space of all $\N_0\cup\{\infty\}$-valued Radon
measures on $X$ endowed with the vague topology.

\begin{th}\label{uigvfuz} The semigoup  $({\bf P}_t)_{t\ge0}$ is
ergodic in the following sense\rom: For an arbitrary probability
measure $\mu$ on $(\Gamma_X,{\cal B}(\Gamma_X))$\rom, the image
measure   ${\bf P}_t\mu$ of $\mu$ under ${\bf P}_t$  converges to
$\pi_m$ weakly on $\dd\Gamma_X$ as $t\to\infty$\rom.
\end{th}

\noindent{\it Proof}. It is evidently enough to prove that, for
each $\gamma\in\Gamma_X$, ${\bf P}_{t,\gamma}$ converges to
$\pi_m$ weakly on $\dd\Gamma_X$ as $t\to\infty$. To this end, it
suffices to show that, for each $\varphi\in C_0(X)$,
$\varphi\le0$,
\begin{equation}\int_{\Gamma_X}e^{\la \varphi,\eta\ra}\,{\bf P}_{t,\gamma}(d\eta)
\to \exp\bigg(\int_X(e^{\varphi(x)}-1)\,m(dx)\bigg)\quad \text{as
$t\to\infty$, for each
$\gamma\in\Gamma_X$},\label{gzszgi}\end{equation}
 see \eqref{ewrwerewrwe5}
and \cite{Ka75}, Chapter~4, in particular, Theorem~4.2. But
analogously to \eqref{hdfvbei}, we get:
\begin{equation}\int_{\Gamma_X}e^{\la \varphi,\eta\ra}\,{\bf
P}_{t,\gamma}(d\eta) =
\prod_{x\in\gamma\cap\operatorname{supp}\varphi}(e^{-t}(e^{\varphi(x)}-1)+1)
\exp\bigg((1-e^{-t})\int_X(e^{\varphi(x)}-1)\,m(dx)\bigg).
 \label{zgusdg}\end{equation}
Since the support of $\varphi$ is compact, we have
$|\gamma\cap\operatorname{supp}\varphi|<\infty$, which implies
that the right hand side of \eqref{zgusdg} converges to the right
hand side of \eqref{gzszgi} as $t\to\infty$.\quad $\square$

\section{Glauber dynamics}\label{ucdguz}

A Markov process on $\Gamma_X$ with the generator $H$ we shall
call a Glauber dynamics of a continuous system of free particles.

\begin{th}\label{tfut6t89} The Glauber dynamics of a continuous system of free particles
 may be realized
as the unique\rom, time homogeneous Markov process $${\bf
M}=({\pmb{ \Omega}},{\bf F},({\bf F}_t)_{t\ge0},({\pmb
\theta}_t)_{t\ge0}, ({\bf P }_\gamma)_{\gamma\in\Gamma_X},({\bf
X}_t)_{t\ge0})$$ on the state space $(\Gamma_X,{\cal
B}(\Gamma_X))$ with transition probability function $({\bf
P}_t)_{t\ge0}$ and with cadlag paths\rom, i\rom.e\rom{.,}  right
continuous on $[0,\infty)$ and having left limits on $(0,\infty)$
\rom(cf\rom.\ e\rom.g\rom.\ \rom{\cite{BG}).}
\end{th}

\noindent{\it Proof}. In what follows, for a metric  space $E$, we
shall   denote by $D([0,\infty),E)$ the space of all cadlag
functions from $[0,\infty)$ into $E$ equipped with the
corresponding cylinder $\sigma$-algebra.

Now, we consider on the two-point space $\{0,1\}$ the time
homogeneous Markov process with cadlag paths, whose transition
probabilities are given by $p_{t,1}(\{1\})=e^{-t}$,
$p_{t,1}(\{0\})=1-e^{-t}$, and $p_{t,0}(\{0\})=1$, $t\ge0$. Let
$\Omega:=D([0,\infty),\{0,1\})$ and let $\bar p$ denote the
probability measure on $\Omega$ defined as the law of the above
process starting at $1$.

For an arbitrary configuration $\gamma\in\Gamma_X$,
$\gamma\ne\varnothing$, on the product space $\Omega^\gamma$ we
consider
 the product measure $\bigotimes _{x\in\gamma}\bar
p_x$, where $\bar p_x:=\bar p$, $x\in\gamma$. We define
\begin{equation}\label{ftf} \Omega^\gamma
\ni\omega=(\omega_x)_{x\in\gamma}\mapsto
Y(\omega)=Y(\omega)(t):=\sum_{x\in\gamma}\omega_x(t)\varepsilon_x,\quad
t\ge0 .\end{equation} As easily seen, for each
$\omega\in\Omega^\gamma$, $Y(\omega)$ is an element of $\pmb
\Omega:=D([0,\infty),\Gamma_X)$ and, furthermore, the mapping
$Y$ is measurable. Let $P_\gamma^Y$ denote the image of
$\bigotimes _{x\in\gamma}\bar p_x$ under $Y$. We also set
$P^Y_\varnothing$ to be the delta measure concentrated at
$Y(t)\equiv\varnothing$, $t\ge0$.

Next, we consider the configuration space $\Gamma_{X\times\R_+}$
over $X\times\R_+$, where $\R_+:=(0,+\infty)$. Let
$\pi_{m\otimes d\tau}$ denote the Poisson measure on ${\cal
B}(\Gamma_{X\times\R_+})$ with intensity measure $m\otimes d\tau$,
where $d\tau$ denotes  the Lebesgue measure on $\R_+$. By
$\widetilde \Gamma_{X\times\R_+}$ we denote the subset of
$\Gamma_{X\times\R_+}$ which consists of those configurations
$\xi$ for which $(x,t),(x',t')\in\xi$, $(x,t)\ne (x',t')$ implies
$x\ne x'$. It is not hard to see that $\widetilde
\Gamma_{X\times\R_+}\in {\cal B}(\Gamma_{X\times\R_+})$.
Furthermore, since the measure $m$ is non-atomic, we have that
$\widetilde \Gamma_{X\times\R_+}$ is of full $ \pi_{m\otimes
d\tau}$ measure.

Now, for any $\xi\in\widetilde\Gamma_{X\times\R_+}$,
$\xi\ne\varnothing$, we consider on the product space $\Omega^
\xi$ the product measure $\bigotimes _{(y,\tau)\in\xi}\bar
p_{(y,\tau)}$, where $\bar p_{(y,\tau)}\equiv \bar p$, and define
\begin{equation}\label{zqweg87e} \Omega^\xi \ni \omega=(\omega_{(y,\tau)})_
{(y,\tau)\in\xi}\mapsto
Z(\omega)=Z(\omega)(t):=\sum_{(y,\tau)\in\xi}\pmb1_{[\tau,\infty)}(t)\omega_{(y,\tau)}(t-\tau)
\varepsilon _y.\end{equation} Again, for each
$\omega\in\Omega^\xi$, $Z(\omega)\in \pmb\Omega$, $Z$ is
measurable, and  by $P_\xi^Z$ we denote the probability measure  on
$\pmb\Omega$ that is the law of the process $Z$. Furthermore, a
monotone class argument shows that, for each measurable set $A$ in
$\pmb\Omega$, the mapping $\widetilde
\Gamma_{X\times\R_+}\ni\xi\mapsto P_\xi^Z(A)\in\R$ is measurable.
Therefore, we can define the probability measure
\begin{equation}\label{bgzug} P^Z:=\int_{\widetilde\Gamma_{X\times\R_+}}\pi_{m\otimes
d\tau}(d\xi)P_\xi^Z\end{equation} on $(\pmb\Omega,{\bf F})$, where
$\bf F$ denotes the cylinder $\sigma$-algebra on $\pmb \Omega $.

Analogously to \eqref{gzgzuguguz}, we define, for each
$\gamma\in\Gamma_X$, ${\bf P_\gamma}$ as the probability measure
on $(\pmb\Omega,{\bf F})$ given by \begin{equation}\label{sukus}
{\bf P}_\gamma (A):=\int_{\pmb\Omega}
P_\gamma^Y(d\omega_1)\int_{\pmb\Omega}P^Z(d\omega_2)\,\pmb1_A(\omega_1+\omega_2),\qquad
A\in{\bf F}.\end{equation} (Notice that
$\omega_1+\omega_2\in\pmb\Omega$ for $P_\gamma^Y\otimes P^Z$-a.e.\
$(\omega_1,\omega_2)\in\pmb\Omega$.) We evidently have that $
\omega(0)=\gamma$ for ${\bf P}_\gamma$-a.e.\
$\omega\in\pmb\Omega$.

Define now ${\bf X}_t(\omega):=\omega(t)$,
$\omega\in\pmb\Omega$, $t\ge0$, and  ${\bf F}_t:=\sigma\{{\bf
X}_s,\, 0\le s\le t\}$. Define the translations $(\pmb\theta
\omega)(s):=\omega(s+t)$, $\omega\in\pmb\Omega$, $s,t\ge0$, and  $${\bf M}:=({\pmb{ \Omega}},{\bf F},({\bf
F}_t)_{t\ge0},({\pmb \theta}_t)_{t\ge0}, ({\bf P
}_\gamma)_{\gamma\in\Gamma_X},({\bf X}_t)_{t\ge0}).$$

To show that $\bf M$ is a time homogeneous Markov process on
$(\Gamma_X,{\cal B}(\Gamma_X))$ with transition probability
function $({\bf P}_t)_{t\ge0}$, it suffices to show that the
finite-dimensional distributions of ${\bf X}_t$ are determined by
the Markov semigroup of kernels $({\bf P}_t)_{t\ge0}$ (see e.g.\
\cite[Ch.~1, Sect.~3]{Bl}).

By \eqref{bgzug} and \eqref{sukus}, we get, for any
$\gamma\in\Gamma_X$, $0< t_1< t_2<\dots<t_n$,
$\varphi_1,\dots,\varphi_n\in{\cal C}$, $n\in\N$,
\begin{align} &\int_{\pmb\Omega}{\bf P}_\gamma(d\omega)
\prod_{k=1}^n
\exp[\la\log(1+\varphi_k),\omega(t_k)\ra]\notag\\&\quad=\int_{\pmb\Omega}
P_\gamma^Y(d\omega_1) \prod_{k=1}^n
\exp[\la\log(1+\varphi_k),\omega_1(t_k)\ra]\notag\\ &\qquad\times
\int_{\widetilde\Gamma_{X\times\R_+}}\pi_{m\otimes
d\tau}(d\xi)\int_{\pmb\Omega}P^Z_\xi(d \omega_2)\prod_{k=1}^n
\exp[\la\log(1+\varphi_k),\omega_2(t_k)\ra].\label{uweizop}
\end{align}
By \eqref{ftf} and the construction of the measure $P_{t,\gamma}$,
we evidently get:
\begin{align}&\int_{\pmb\Omega}
P_\gamma^Y(d\omega_1) \prod_{k=1}^n
\exp[\la\log(1+\varphi_k),\omega_1(t_k)\ra]\notag\\ &\qquad=
\int_{\Gamma_X}P_{t_1,\gamma}(d\eta_1)\exp[\la\log(1+\fii_1),\eta_1\ra]
\int_{\Gamma_X}P_{t_2-t_1,\eta_1}(d\eta_2)\exp[\la\log(1+\fii_2),\eta_2\ra]\notag\\
&\qquad\quad \times \dotsm\times
\int_{\Gamma_X}P_{t_n-t_{n-1},\eta_{n-1}}(d\eta_n)\exp[\la\log(1+\fii_n),\eta_n\ra].
\label{wetukfcwef}\end{align} For any
$\Delta_1,\dots,\Delta_k\in{\cal B}(X\times\R_+)$, $k\in\N$, which
are of finite $m\otimes d\tau$ measure and disjoint, the random
variables $\la\pmb 1_{\Delta_1},\cdot\ra,\dots,\la\pmb
1_{\Delta_k},\cdot\ra$ are independent under $\pi_{m\otimes
d\tau}$. Therefore, by \eqref{zqweg87e},
\begin{align}&\int_{\widetilde\Gamma_{X\times\R_+}}\pi_{m\otimes
d\tau}(d\xi)\int_{\pmb\Omega}P_\xi^Z(d \omega_2)\prod_{k=1}^n
\exp[\la\log(1+\varphi_k),\omega_2(t_k)\ra]\notag\\&\qquad=\int_{\widetilde\Gamma_{X\times\R_+}}\pi_{m\otimes
d\tau}(d\xi)\prod_{(y,\tau)\in\xi:\,
y\in\bigcup_{k=1}^n\operatorname{supp}\varphi_k} \int_{\Omega}\bar
p(d\omega_{(y,\tau)})\notag\\ &\qquad\quad\times\prod_{k=1}^n
\exp[\la \log(1+\varphi_k),\pmb
1_{[\tau,\infty)}(t_k)\omega_{(y,\tau)}(t_k-\tau)\varepsilon_y\ra]
\notag\\&\qquad=\prod_{k=1}^n\int_{\widetilde\Gamma_{X\times
(t_{k-1},t_k)}}\pi_{m\otimes d\tau}(d\xi)\prod_{(y,\tau)\in\xi}
\int_{\Omega}\bar p(d\omega_{(y,\tau)})
\notag\\&\qquad\quad\times\prod_{l=k}^n \exp[\la
\log(1+\varphi_l),
\omega_{(y,\tau)}(t_l-\tau)\varepsilon_y\ra],\label{hsdjisdui}
\end{align}
where $t_{0}:=0$. Using the Markov property, we get, for any
fixed $k\in\{1,\dots,n\}$ and $(y,\tau)\in X\times (t_{k-1},t_k)$:
\begin{align}
&\int_\Omega \bar p(d\omega) \prod_{l=k}^n
\exp[\la\log(1+\varphi_l),\omega(t_l-\tau)\varepsilon_y\ra]
\notag\\ &\qquad =e^{-(t_k-\tau)}(1+\varphi_k(y))
\int_{\{0,1\}}p_{t_{k+1}-t_k,1}(d
a_{k+1})\exp[\la\log(1+\varphi_{k+1}),a_{k+1}\varepsilon_y\ra]
\notag\\ &\qquad\quad\times\dots\times
\int_{\{0,1\}}p_{t_n-t_{n-1},a_{n-1}}(da_n)
\exp[\la\log(1+\fii_n),a_{n}\varepsilon_y\ra]+(1-e^{-(t_k-\tau)}).
\label{gid}\end{align} Hence, by \eqref{ewrwerewrwe5},
\eqref{hsdjisdui}, and \eqref{gid}, we get:
\begin{align}&\int_{\widetilde\Gamma_{X\times\R_+}}\pi_{m\otimes
d\tau}(d\xi)\int_{\pmb\Omega}P^Z_\xi(d \omega_2)\prod_{k=1}^n
\exp[\la\log(1+\varphi_k),\omega_2(t_k)\ra]\notag\\&\qquad=\exp\bigg[
\sum_{k=1}^n(1-e^ {-(t_k-t_{k-1})})\int_X
m(dy)\,\bigg(-1+(1+\varphi_k(y))\notag\\ &\qquad\quad\times
\int_{\{0,1\}}p_{t_{k+1}-t_k,1}(da_{k+1})\exp[\la\log(1+\varphi_{k+1}),a_{k+1}\varepsilon_y\ra]
\notag\\ &\qquad\quad\times\dots\times
\int_{\{0,1\}}p_{t_n-t_{n-1},a_{n-1}}(da_n)
\exp[\la\log(1+\fii_n),a_{n}\varepsilon_y\ra]
\bigg)\bigg].\label{gdu834u}
\end{align}

It is easy to see that, for any $t>0$ and any
$\gamma_1,\gamma_2\in\Gamma_X$, $\gamma_1\cap\gamma_2=\varnothing
$, the measure $P_{t,\gamma_1+\gamma_2}$ on $\Gamma_X$ is the
convolution of the measures $P_{t,\gamma_1}$ and $P_{t,\gamma_2}$.
Therefore, by \eqref{ewrwerewrwe5}, \eqref{gzgzuguguz},
\eqref{uweizop}, \eqref{wetukfcwef}, and \eqref{gdu834u},
\begin{align}& \int_{\Gamma_X}{\bf
P}_{t_1,\gamma}(d\gamma_1)\exp[\la\log(1+\fii_1),\gamma_1\ra]\int_{\Gamma_X}{\bf
P
}_{t_2-t_1,\gamma_1}(d\gamma_2)\exp[\la\log(1+\fii_2),\gamma_2\ra]\\
&\qquad\quad\times \dots\times \int_{\Gamma_X}{\bf P
}_{t_n-t_1,\gamma_{n-1}}(d\gamma_n)\exp[\la\log(1+\fii_n),\gamma_n\ra]\notag\\
&\qquad
=\int_{\Gamma_X}P_{t_1,\gamma}(d\eta_1)\int_{\Gamma_X}\pi_{(1-e^{t})m}(d\theta_1)\exp[\la\log(1+\varphi_1),
\eta_1+\theta_1\ra]\notag\\ &\qquad\quad \times
\int_{\Gamma_X}P_{t_2-t_1,\eta_1+\theta_1}(d\eta_2)\int_{\Gamma_X}\pi_{(1-e^{-(t_2-t_1)})m}(d\theta_2)
\exp[\la\log(1+\varphi_2),\eta_2+\theta_2\ra]\notag\\
&\qquad\quad\times\dots \times
\int_{\Gamma_X}P_{t_n-t_{n-1},\eta_{n-1}+\theta_{n-1}}(d\eta_n)\notag\\
&\qquad\quad\times\int_{\Gamma_X}\pi_{(1-e^{-(t_n-t_{n-1})})m}(d\theta_n)
\exp[\la\log(1+\varphi_n),\eta_n+\theta_n\ra]\notag\\ &\qquad =
\int_{\Gamma_X}P_{t_1,\gamma}(d\eta_1)\exp[\la\log(1+\varphi_1),\eta_1\ra]\int_{\Gamma_X}
P_{t_2-t_1,\eta_1}(d\eta_2)\exp[\la\log(1+\varphi_2),\eta_2\ra]\notag\\&\qquad\quad
\times\dots\times
\int_{\Gamma_X}P_{t_n-t_{n-1},\eta_{n-1}}(d\eta_{n})\exp[\la\log(1+\varphi_n),\eta_n\ra]\notag\\
&\qquad\quad\times\prod_{k=1}^n\int_{ \Gamma_X}
\pi_{(1-e^{-(t_k-t_{k-1})})m}(d\theta_k)\exp[\la\log(1+\fii_k),\theta_k\ra]\notag\\
&\qquad\quad\times\int_{\Gamma_X}P_{t_{k+1}-t_k,\theta_k}(d\eta_{k+1})\exp[\la\log(1+\fii_{k+1}),\eta_{k+1}\ra]
\notag\\&\qquad\quad\times\dots\times
\int_{\Gamma_X}P_{t_n-t_{n-1},\eta_{n-1}}(d\eta_n)\exp[\la\log(1+\varphi_n),\eta_n\ra]\notag\\
&\qquad=\int_{\pmb\Omega} P_\gamma^Y(d\omega_1) \prod_{k=1}^n
\exp[\la\log(1+\varphi_k),\omega_1(t_k)\ra]\notag\\&\qquad\quad\times\prod_{k=1}^n\int_{
\Gamma_X}
\pi_{(1-e^{-(t_k-t_{k-1})})m}(d\theta_k)\exp[\la\log(1+\fii_k),\theta_k\ra]\notag\\&\qquad\quad
\times
\prod_{x\in\theta_k}(1+\fii_k(x))\int_{\{0,1\}}p_{t_{k+1}-t_k,1}(da_{k+1})\exp[\la\log(1+\fii_{k+1}),
a_{k+1}\varepsilon_x\ra]\notag\\&\qquad\quad\times\dots\times\int_{\{0,1\}}p_{t_n-t_{n-1},a_{n-1}}
(da_n)\exp[\la\log(1+\fii_n),a_n\varepsilon_x\ra]\notag\\
&\qquad=\int_{\pmb\Omega}{\bf P}_\gamma(d\omega) \prod_{k=1}^n
\exp[\la\log(1+\varphi_k),\omega(t_k)\ra].\label{wewewe}
\end{align}

Hence, analogously to the proof of \eqref{uacad}, we conclude from
 \eqref{wewewe} that, for any  $A_1,\dots,A_n\in{\cal B}(\Gamma_X)$, $$
 {\bf P}_\gamma({\bf X}_{t_1}\in A_1,\dots,{\bf X}_n\in
 A_n)=\int_{A_1}{\bf P}_{t_1,\gamma}(d\gamma_1)\int_{A_2}{\bf
 P}_{t,\gamma_1}(d\gamma_2)\dotsm\int_{A_n}{\bf
 P}_{t_n,\gamma_{n-1}}(d\gamma_n).$$

Finally, we note that any measure on the space $({\pmb
\Omega},{\bf F})$ is uniquely  determined by its
finite-dimensional distributions, and therefore the constructed
Markov process with cadlag paths is unique. \quad$\square$

\begin{cor}\label{huau} The Markov process $\bf M$ from
Theorem~\rom{\ref{tfut6t89} } is a strong Markov process\rom.
\end{cor}

\noindent {\it Proof}. The statement follows directly from
Theorems~\ref{haajsd}, \ref{tfut6t89}  and
\cite[Theorem~5.10]{dynkin}.\quad $\square$

\begin{rem}\rom{ It is easy to see that the process $({\bf
X}_t)_{t\ge0}$ constructed in the course of the proof of
Theorem~\ref{tfut6t89}   is even Markov with respect to the
filtration $({\bf F}_{t+})_{t\ge0}$, where ${\bf
F}_{t+}:=\bigcap_{s>t}{\bf F}_t$.}\end{rem}

Suppose that the space $X$ is non-compact and  by
$\Gamma_{X,\,{\mathrm inf}}$ denote the subset  of $\Gamma_X$ consisting
of all infinite configurations $\gamma\in\Gamma_X$. It is not hard
to show that $\Gamma_{X,\,{\mathrm inf}}\in{\cal B}(\Gamma_X)$. We
endow $\Gamma_{X,\,{\mathrm inf}}$ with the vague topology and
denote by ${\cal B}(\Gamma_{X,\,{\mathrm inf}})$ the corresponding
Borel $\sigma$-algebra.

\begin{cor} \label{wfdhci} Suppose that $m(X)=\infty$\rom. Then\rom, in the formulation of
Proposition~\rom{\ref{cjj2356},} Theorems~\rom{\ref{chbahs},
\ref{haajsd}, \ref{tfut6t89},}
 and Corollary~\rom{\ref{huau},}
we can replace $\Gamma_X$ by  $\Gamma_{X,\,{\mathrm inf}}$\rom.
\end{cor}

\noindent {\it Proof}. It is well known that the condition
$$m(X)=\infty$$ implies that $\pi_{zm}(\Gamma_{X,\,{\mathrm
inf}})=1$, $z>0$. Furthermore, by \eqref{gzgzuguguz}, $$ {\bf
P}_{t,\gamma}(\Gamma_{X,\,{\mathrm inf}})=1,\qquad t>0,\ \gamma\in
\Gamma_{X,\,{\mathrm inf}}.$$ Therefore, by
Proposition~\ref{cjj2356}, $({\bf P}_t)_{t\ge 0}$ is a Markov
semigroup of kernels on $(\Gamma_{X,\,{\mathrm inf}},{\cal
B}(\Gamma_{X,\,{\mathrm inf}}))$, and we can replace $\Gamma_X$
by  $\Gamma_{X,\,{\mathrm inf}}$ in Theorems~\ref{chbahs},
\ref{haajsd}.

Next, by using the Borel-Cantelli lemma, we easily get:
\begin{equation}\label{jnasj}P_{t,\gamma}(\Gamma_{X,\,{\mathrm
inf}})=1,\qquad t>0,\ \gamma\in \Gamma_{X,\,{\mathrm inf}}.\end
{equation} By the construction of $P^Y_\gamma$, we also have:
\begin{equation}\label{dtrd} P^Y_\gamma(\omega\in\pmb\Omega:\,
\omega(\tau)\supset \omega(t)\,\forall \tau\in[0,t])=1,\qquad
t\in\N,\ \gamma\in \Gamma_{X,\,{\mathrm inf}}.\end{equation}

Taking \eqref{jnasj} by  \eqref{dtrd} into account, we can now
replace $\Gamma_X$ by  $\Gamma_{X,\,{\mathrm inf}}$ in the proof
of Theorem~\ref{tfut6t89}, to get the respective modification of
the latter theorem. \quad $\square$

\begin{rem}\rom{ Notice that $D([0,\infty),\Gamma_{X,\,{\mathrm
inf}})$ does not belong to the cylinder $\sigma$-algebra
on\linebreak $D([0,\infty),\Gamma_X)$, so that in the proof of
Corollary~\ref{wfdhci} we could not simply state that\linebreak
$D([0,\infty),\Gamma_{X,\,{\mathrm inf}})$ is a set of full ${\bf
P }_{\gamma}$ measure for each $\gamma\in\Gamma_{X,\,{\mathrm
inf}}$.

}\end{rem}


\begin{thebibliography}{99}

\bibitem{AKR3}  S.  Alberverio, Yu.\ G. Kondratiev, and M. R\"ockner, Analysis and
  geometry on configuration spaces, {\it J. Func.\ Anal.}\ {\bf 154} (1998), 444--500.

\bibitem{BK} Yu.\ M. Berezansky and  Yu.\ G. Kondratiev, ``Spectral Methods in
  Infinite Dimensional Analysis,''  Kluwer Acad.\ Publ.,  Dordrecht, 1994.

\bibitem{BCC} L. Bertini, N. Cancrini, and F. Cesi,
The spectral gap for a Glauber-type dynamics in a continuous gas,
{\it  Ann.\ Inst.\ H. Poincar\'e Probab.\ Statist.}\  {\bf  38} (2002),
91--108.


\bibitem{Bl} R. M. Blumenthal, ``Excursions of Markov Processes,''
Birkh\"auser, Boston, 1992.

\bibitem{BG} R. M. Blumenthal and R. K. Getoor, ``Markov Processes
and Potential Theory,'' Academic Press, New York, 1968.


\bibitem{dynkin} E. B. Dynkin, ``Theory of Markov Processes,''
Pergamon Press, Oxford,  1960.


\bibitem{HS} R. A. Holley and  D. W.  Stroock,
 Nearest neighbor birth and death processes on the real line, {\it Acta Math.}\ {
\bf 140} (1987), 103--154.

\bibitem{Ka75} O. Kallenberg, ``Random Measures,'' Academic Press,
San Diego, 1975.


\bibitem{KMM}  J. Kerstan,  K.~Matthes, and J.~Mecke, ``Infinite
Divisible Point Processes,'' Akademie-Verlag, Berlin, 1978.

\bibitem{Kingman} J. F. C.  Kingman,  ``Poisson processes,'' Oxford University Press,
Oxford, 1993.


\bibitem{KL} Yu.~Kondratiev and E.~Lytvynov,   Glauber dynamics of continuous particle
systems, to appear in {\it  Ann.\ Inst.\ H. Poincar\'e Probab.\ Statist.}


\bibitem{KLR}  Yu.~Kondratiev, E. Lytvynov,
and  M.~R\"ockner, The heat semigroup on configuration spaces,
{\it Publ.\ Res.\ Inst.\ Math.\ Sci.} {\bf 39} (2002), 1--48.


\bibitem{MR1}  Z.-M. Ma and M.   R\"ockner, ``An Introduction to the
Theory of (Non-Symmetric) Dirichlet Forms,'' Springer-Verlag, Berlin, 1992.

\bibitem{MR2}  Z.-M. Ma and M.   R\"ockner, Construction of
diffusions on configuration spaces,  {\it Osaka J. Math.}\  {\bf 37} (2000),
273--314.



\bibitem{Me67} J. Mecke, Station\"are zuf\"allige Ma\ss e auf lokalkompakten Abelschen Gruppen,
{\it Z.
   Wahrscheinlichkeitstheorie und Verw. Gebiete} {\bf 9} (1967),
   36--58.

\bibitem{P} C. Preston,  Spatial birth-and-death processes.  With discussion, {\it in}
 ``Proceedings of the 40th Session of the International Statistical
Institute (Warsaw, 1975), Vol.~2,''  Bull.\ Inst.\ Internat.\
Statist., Vol.~46, 1975, pp.~371--391.


\bibitem{RS}  M. Reed and B. Simon,  ``Methods of Modern Mathematical
Physics, Vol.~1. Functional Analysis,''  Academic Press, San Diego,
1972.

\bibitem{S1} D. Surgailis, On multiple Poisson stochastic integrals and associated Markov semigroups,
{\it  Probab.\ Math.\
   Statist.}\ {\bf 3} (1984), 217--239.

\bibitem{S2} D. Surgailis, On Poisson multiple stochastic integrals and associated equilibrium Markov processes,
{\it in}   ``Theory
   and Application of Random Fields (Bangalore, 1982),''  Lecture Notes in
   Control and Inform. Sci., Vol.~49, Springer, Berlin,
  1983, pp.~233--248.





\bibitem{Wu} L. Wu, Estimate of spectral gap for continuous gas,
Preprint, 2003.



  \end{thebibliography}
\end{document}